\documentclass[12pt]{article}

\usepackage{amsmath}
\usepackage{amsfonts}
\usepackage{amsthm}

\numberwithin{equation}{section}

\def\Cmpx{{\mathbb{C}}}
\def\Real{{\mathbb{R}}}
\def\Intg{{\mathbb{Z}}}
\def\cnj{\overline}

\def\inter{{\cap}}

\def\tfrac#1#2{{\textstyle{\frac{#1}{#2}}}}
\def\half{{\textstyle{\frac{1}{2}}}}
\def\nfrac#1#2{{\raisebox{.5ex}{$#1$}\!/\!\raisebox{-.5ex}{$#2$}}}
\def\scrhalf{{\raisebox{.3ex}%
{$\scriptstyle 1$}\!/\!\raisebox{-.3ex}{$\scriptstyle 2$}}}

\def\wt{\widetilde}
\def\wh{\widehat}

\def\spanrm{{\mathrm{span}}}
\def\innerprod(#1){{\langle #1 \rangle}}
\def\norm#1{\|#1\|}
\def\rvec(#1){{|#1\rangle}} 
\def\lvec(#1){{\langle#1|}}

\def\defn#1{{\em #1}}

\newtheorem{theorem}{Theorem}
\newtheorem{lemma}[theorem]{Lemma}
\newtheorem{corol}[theorem]{Corollary}

\def\AWeil{{{\cal W}_A}}
\def\BWeil{{{\cal W}_B}}
\def\Weil{{{\cal W}}}
\def\spinset{{{\cal S}}}
\def\PsiSpace{{\Psi}}
\def\PPexpr{{\cal P}}
\def\Hilb{{\cal H}}
\def\baseterm(#1,#2,#3){{\psi{}^{(#2,#3)}_{#1}}}
\def\PsinSpace#1{{{\Psi}{}^{(\bullet,\bullet)}_{#1}}}
\def\PsirSpace#1{{{\Psi}{}^{(#1,\bullet)}_{\bullet}}}
\def\PsimSpace#1{{{\Psi}{}^{(\bullet,#1)}_{\bullet}}}
\def\PsinrSpace#1#2{{{\Psi}{}^{(#2,\bullet)}_{#1}}}
\def\Rp{{\hat R}}
\def\Ad#1{{\mathrm{Ad}}_{#1}}
\def\tr{{\mathrm{tr}}}

\def\apm{a_\pm}
\def\bpm{b_\pm}
\def\wmem{{w}}

\title{A Natural Basis for Spinor and Vector Fields on the
Noncommutative sphere}

\author{{\bf J. Gratus%
\thanks{Funded by the Royal Society of London European 
Science Exchange Programme}}%
\\
Laboratoire de Gravitation et Cosmologie Relativistes%
\thanks{\it Laboratoire associ\'e au CNRS {\rm URA 769}}
\\
Tour 22/12 4eme etage, Boite Courrier142, 4pl Jussieu. F75252
Paris
\\
email: gratus@ccr.jussieu.fr
}

\begin{document}

\maketitle

\begin{abstract}

The product of two Heisenberg-Weil algebras contains the
Jordan-Schwinger representation of $su(2)$. This Algebra is quotiented
by the square-root of the Casimir to produce a non-associative algebra
denoted by $\Psi$. This algebra may be viewed as the right-module over
one of its associative subalgebras which corresponds to the algebra of
scalar fields on the noncommutative sphere. It is now possible to
interpret other subspaces as the space of spinor or vector fields on
the noncommutative sphere. A natural basis of $\Psi$ is given which
may be interpreted as the deformed entries in the rotation matrices of
$SU(2)$.

\end{abstract}


\tableofcontents

\section{Introduction}

The noncommutative or ``fuzzy'' sphere has been considered by several
authors in different contexts, including the general
quantisation procedure, coherent states, noncommutative geometry, 
the theory of membranes, and the quantum Hall effect. 
(See references in
\cite{Gratus5,Madore_book,Watamura1,Grosse5})

Normally the approximation for the algebra of functions on a sphere
(scalar fields) is in terms of matrices,
\cite{Madore_book,Madore2,Grosse4} where the functions on a sphere
appear only in the limit as the size of the matrix tends to infinity.
An alternative approach was presented in \cite{Gratus5} as the a two
parameter algebra of polynomials $\PPexpr(\varepsilon,R)$ with
$\varepsilon,R\in\Real$.

The standard way of defining a vector field is to consider it as a
derivation on the algebra of functions. If the latter is replaced by a
matrix algebra then all such derivations are the adjoint action of
elements in that algebra \cite{Madore_book,Gratus4}.  If
$\{x^1,x^2,x^3\}$ are the coordinate functions for the noncommutative
sphere then the corresponding derivations $X_i=\Ad{x^i}$ obey the
equation $\sum_{i=1}^3 x^i X_i=0$ only in the commutative limit.  Here
we give an alternative definition for the analogue of vector fields
which solves this equation identically in the noncommutative case, but
these vector fields are derivations only in the limit $\varepsilon=0$.

The idea of using the Jordan-Schwinger representation of $su(2)$ to
describe spinor fields have be pursued by many authors
\cite{Watamura1,Grosse5,Grosse1,Grosse3,Grosse6}. 
They discovered that one could view spinor fields as rectangular
matrices. 

This article may be consider as an extension of the \cite{Gratus5} to
spinor and vector fields by the use of the Jordan-Schwinger
representation of $su(2)$.


\subsection{Structure of the article}

We present $\Weil$ an algebra which is the
product of two Heisenberg-Weil algebra generated by
$[a_-,a_+]=\varepsilon$ and $[b_-,b_+]=\varepsilon$.  This naturally
contains the Jordan-Schwinger representation of both $su(2)$ given
by $\{J_0,J_+,J_-\}$ and $su(1,1)$ given by $\{K_0,K_+,K_-\}$ defined
in (\ref{W_def_JK}). 

The idea is to consider the algebra generated by $\{J_0,J_+,J_-\}$ to
be equivalent to the analogue of the algebra of functions (or scalar
fields) on the sphere as presented in \cite{Gratus5}. In this paper
the algebra generated by $\{J_0,J_+,J_-\}$ which is the universal
covering algebra of $su(2)$ is quotiented by the ideal generated the
Casimir $J_0^2 + \tfrac12 J_+J_- + \tfrac12 J_-J_+ \sim R^2$.  where
$R^2\in\Real$. Since the Casimir operator commutes with the set
$\{J_0,J_+,J_-\}$ this ideal is two sided and the corresponding
quotient is an algebra.

In this paper by contrast the Casimir operator is given by $J_0^2 +
\tfrac12 J_+J_- + \tfrac12 J_-J_+ = K_0^2 - \tfrac14 \varepsilon^2$.
We would like to take the square root of this equation and consider
quotienting $\Weil$ by the ideal generated by
$K_0\sim\Rp=(R^2+\tfrac14\varepsilon^2)^\scrhalf$. However although
$K_0$ commutes with any polynomial in $\{J_0,J_+,J_-\}$, it does not
commute with all the elements in $\Weil$. Therefore the left ideal
generated by $K_0\sim\Rp$ is not two sided and the corresponding
quotient product is not associative.  However our persistence is
rewarded.

We define the subspace $\PsiSpace\subset\Weil$ as all the totally
symmetric polynomials in $(a_+,a_-,b_+,b_-)$, 
which contain no factor of $K_0$. This is achieved by defining the
\defn{formal trace}.  We can define two products on $\PsiSpace$,
$\rho$ and $\rho^\star$ by quotienting by either the left or right
ideal generated by $K_0\sim\Rp$. 

There is a natural basis of $\PsiSpace$ given by the set
$\{\baseterm(n,r,m)\ |\ 2n,n+r,n+m\in\Intg, n\ge0, |r|\le n |m|\le
n\}$. The basis element $\baseterm(n,r,m)$ is the unique (up to scale
factor) homogeneous polynomial of degree $2n$, which is an eigenstate
of both $\Ad{K_0}$ and $\Ad{J_0}$ with corresponding eigenvalues
$\varepsilon r$ and $\varepsilon m$. These basis elements are also
orthogonal with respect to a sesquilinear form defined in a similar
way to the sesquilinear form in \cite{Gratus5}.

These basis elements may be thought of as entries in a deformed
rotation matrix, since, in the limit $\varepsilon=0$,
$\baseterm(n,r,m)$ is proportional to $D^n_{mr}(\alpha,\beta,\gamma)$
where $\alpha,\beta,\gamma$ are the Eular angles. As a result we may
interpret the algebra $(\PsiSpace,\rho)$ as the noncommutative and
nonassociative analogue of functions on the group $SU(2)$.

We label the eigenspace
$\PsirSpace{r}=\spanrm\{\baseterm(n,r,m)\forall n,m\}$.  Since the
eigenspace $\PsirSpace{0}$ corresponds to all elements of $\PsiSpace$
which commute with $K_0$, these elements can be written as polynomials
in $\{J_0,J_+,J_-\}$ and thus correspond to the scalar fields on the
noncommutative sphere.  Now as stated the product $\rho$ on
$\PsiSpace$ is not associative but we show that $\PsirSpace{r}$ for
each $r$ can be viewed as a right module over $\PsirSpace{0}$.  This
is analogous to the fact that standard spinor and vector fields are
modules over the space of scalar fields.

It is natural to call the set $\PsirSpace{1}$ the analogue of vectors
fields, and to call $\PsirSpace{-1}$ the analogue of covectors fields
or 1-forms.  Using this interpretation we may define an exterior
``derivative'' $d:\PsirSpace{0}\mapsto\PsirSpace{-1}$. This map
obeys the Leibniz rule only in the limit $\varepsilon=0$. Likewise the
interpretation of vectors as derivations on the algebra of function is
also true only in this limit. By contrast the equation $x^iX_i=0$ is
true for all $\varepsilon$ not just in the limit.  Thus its swings and
roundabouts and depends on ones personal conviction as to which
properties of a vector are fundermental and thus should be valid for
all $\varepsilon$ and which need be true only in the limit.

The spaces $\PsirSpace{1/2}$ and $\PsirSpace{1/2}$, are the
analogue of spinor fields on the noncommutative sphere. They may be
viewed as contra-variant and co-variant vector fields, or positive and
negative chiral vector fields depending on ones interpretation.  We
justify the definition of a spinor field on two grounds. First the
product of two spinor fields is a vector, and secondly the rotation of
a spinor though $2\pi$ inverses the sign of the spinor.
Alternatively writing the spinors as 2-dimensional vectors with
entries in $\PsirSpace{-1/2}\oplus\PsirSpace{1/2}$, we obtain the
standard definition of spinors in the limit $\varepsilon=0$.

\subsection{Contents of the article}

Most of this article concerns itself with the mathematical structure
necessary to define $\PsiSpace$ and its products.  In section
\ref{ch_Weil} we examine the Heisenberg-Weil algebra $\Weil$.  We
define the formal trace and give some of its properties.  In section
\ref{ch_Pex} we present the vector space $\PsiSpace$ and some of its
subspaces. We define the non associative products $\rho$ and
$\rho^\star$, and show that when they are restricted to the space of
scalar fields, they are associative.  We also define the sesquilinear
Hermitian product that fails to be positive definite. (Although it is
positive definite in the case $\varepsilon=0$.)

In section \ref{ch_Pnrm} we define the orthogonal basis for
$\PsiSpace$ in a similar vain to that in \cite{Gratus5}. 
In section \ref{ch_form} we give useful formulae for
calculating the product $\rho$ on $\PsiSpace$. We show how to write
the basis elements either as rectangular matrices or in terms of 
Hahn polynomials.

In section \ref{ch_Phys} we demonstrate how $\baseterm(n,r,m)$ may be
regarded as the deformed rotation matrix entry, and also how
$\PsirSpace{1}$ can be viewed as the space of vector fields and
$\PsirSpace{-1}$ the space of covector fields, and $\PsirSpace{1/2}$
as the space of spinor fields on the noncommutative sphere.

Finally in section \ref{ch_Prob} we discuss some of the problems with
this interpretation. We also suggest how this methodology can be
extended to scalar, spinor and vector fields on more exotic symmetric
spaces such as the Einstein-DeSitter universe.

\subsection{Notations}

The summation convention is not used in this article.

The notational difference with \cite{Gratus5} are that (1) the basis
elements are not normalised, (2) we use $\varepsilon$ instead of
$\kappa$, and (3) since the product is not associative
it is written explicitly whenever used.

$\Rp$ is always given in relation to $R$ by (\ref{Pex_def_Rp}).  In
most cases we do not write the dependence on $R$ and $\varepsilon$
explicitly unless there is room for confusion. The implicit
dependencies on $\Rp$ and $\varepsilon$ are given by:

\begin{tabular}{p{10cm}cc}
& $\Rp$ & $\varepsilon$
\\
$\Weil$ as a vector space & No & No \\
$\Weil$ as an algebra & No & Yes \\
$S(s,t,u,v)$ & No & No \\
$\PsiSpace$, $\PsinSpace{n}$, $\PsirSpace{r}$, $\PsimSpace{m}$ as
vector spaces & No & No \\
$\rho:\Weil\mapsto\PsiSpace$ as a projection & Yes & No \\
$\rho:\PsiSpace\times\PsiSpace$ as a product & Yes & Yes \\
$\baseterm(n,r,m)$ & No & No
\end{tabular}

\noindent
We shall use the following  notation for elements of each set:

\begin{tabular}{cp{8cm}cc}
Space & Description & General elements & Basis 
\\
$\Weil$ & Polynomials in $(a_+,a_-,b_+,b_-)$ &
$\wmem,\wmem_1,\wmem_2,\ldots$ & 
$S(s,t,u,v)$ 
\\
$\PsiSpace$ & formally traceless symmetric polynomials &
$\xi,\zeta,\ldots$ &
$\baseterm(n,r,m)$
\\
$\PsirSpace{0}$ & analogue of functions &
$f,g,\ldots$ & 
$\baseterm(n,0,m)$
\\
$\PsirSpace{1}$ & analogue of vectors fields &
$X,Y,\ldots$ &
$\baseterm(n,1,m)$
\\
$\PsirSpace{-1}$ & analogue of 1-forms &
$\xi,\zeta,\ldots$ &
$\baseterm(n,-1,m)$
\\
$\PsirSpace{1/2}$ & analogue of spinors &
$\xi,\zeta,\ldots$ &
$\baseterm(n,1/2,m)$
\end{tabular}
\section{$\Weil$-expressions: The Free Algebra of the Product of Two
Heisenberg-Weil algebras}
\label{ch_Weil}

As stated in the introduction in this section we define the algebra
$\Weil=\Weil(\varepsilon)$ which is the product of two Heisenberg-Weil
algebras. This naturally contains the Jordan Schwinger representation
of $su(2)$ and $su(1,1)$.  We define a basis for this algebra in terms
totally symmetric polynomials.  Since we wish quotient out by
$K_0=\Rp$ we need to find representative of the $\Weil$ which contain
no factor of $K_0$. This is achieved by defining the \defn{formal
trace}.  We give some basic properties of the formal trace including
that it commutes with the adjoint representation of $su(2)$. Finally
we show a convenient way of writing the elements of $\Weil$ as non
symmetric polynomials

\vspace{1 em}

For $\varepsilon\in\Real$ let
\begin{align}
\Weil(\varepsilon)=\{\mbox{free algebra of finite polynomials in
$(a_+,a_-,b_+,b_-)$}\}\Big/\sim
\end{align}
where
\begin{align}
[a_-,a_+] &\sim \varepsilon &
[b_-,b_+] &\sim \varepsilon &
[a_\pm,b_\pm] &\sim 0 
\label{W_comm_rel}
\end{align}
This may be viewed as the tensor product of two copies of the
Heisenberg-Weil algebra
$\Weil(\varepsilon)=\AWeil(\varepsilon)\otimes\BWeil(\varepsilon)$
where $\AWeil(\varepsilon)$ is the algebra of polynomials of
$(a_+,a_-)$, and likewise for $\BWeil(\varepsilon)$.

{}From now on we write $\Weil=\Weil(\varepsilon)$ when there is no doubt
about $\varepsilon$.  A natural basis for $\Weil$ is given by
$S(s,t,u,v)$ which is the totally symmetric homogeneous polynomial
with each term a permutation of $a_-^s\,a_+^t\,b_-^u\,b_+^v$ and
having coefficient 1.

\vspace{1em}

Let $\Weil^n$ for $2n\in\Intg$, $n\ge0$ be the space of all symmetric
homogeneous polynomials of degree $2n$. The basis for $\Weil^n$ is now
given by $\{S(s,t,u,v)\ |\ s+t+u+v=2n\}$. It is easy to show the
dimension of $\Weil^n$ is given by
$\dim(\Weil^n)=\tfrac16(2n+1)(2n+2)(2n+3)$.

The subspace $\Weil^1$ has dimension 10. Out of these 6 elements turn
out to be very important and are given special names: (written in
non-symmetric by a more convenient form)
\begin{equation}
\begin{aligned}
J_0 &= \tfrac12(a_+a_- - b_+b_-) & 
J_+ &= a_+b_- & 
J_- &= a_-b_+ \\
K_0 &= \tfrac12(a_+a_- + b_+b_- + \varepsilon) & \qquad
K_+ &= a_+b_+ & \qquad
K_- &= a_-b_- 
\end{aligned}
\label{W_def_JK}
\end{equation}
It is easy to show that $\{J_0,J_+,J_-\}$ forms a representation of
$su(2)$ with $[J_0,J_+]=\varepsilon J_+$ etc, and that
$\{K_0,K_+,K_-\}$ form a representation of $su(1,1)$ with
$[K_0,K_+]=\varepsilon K_+$. The Casimir of these are given by
\begin{align}
J_0^2 + \tfrac12 J_+J_- + \tfrac12 J_-J_+ &= K_0^2 - 
\tfrac14 \varepsilon^2 
\label{W_Casimir_J}
\\
K_0^2 - \tfrac12 K_+K_- - \tfrac12 K_-K_+ &= J_0^2 - 
\tfrac14 \varepsilon^2 
\label{W_Casimir_K}
\end{align}
Furthermore $K_0$ commutes with $\{J_0,J_+,J_-\}$ and $J_0$ commutes
with $\{K_0,K_+,K_-\}$. The entire set $\Weil^1$ forms
a representation of $so(3,2)$ \cite{Kibler1}.

\vspace{1em}

There is a \defn{conjugation} of elements in $\Weil$ given by
\begin{align}
\dagger &: \Weil \mapsto \Weil \cr
(\wmem_1 \wmem_2)^\dagger &= \wmem_2^\dagger \wmem_1^\dagger \qquad
(a_\pm)^\dagger = a_\mp \qquad
(b_\pm)^\dagger = b_\mp \qquad
\lambda^\dagger = \cnj{\lambda} \qquad
\forall \wmem_1,\wmem_2\in\Weil,\lambda\in\Cmpx
\end{align}
On the basis elements it is easy to show that
\begin{align}
S(s,t,u,v)^\dagger = S(t,s,v,u)
\end{align}

\vspace{1em}

We define the \defn{formal trace} of an element $\wmem\in\Weil$ as
follows:
\begin{align}
\tr:\Weil\mapsto\Weil,\quad:\Weil^n\mapsto\Weil^{n-1}
\end{align}
Write $\wmem$ as a totally symmetric polynomial. There are sixteen
possible combinations for the last two elements of each term. 
Collecting these terms together, we can now write
\begin{align*}
\wmem = 
\wmem_1 a_+a_- + \wmem_2 a_-a_+ + \wmem_3 b_+b_- + \wmem_4 b_-b_+ +
\wmem_5 a_+^2 + \cdots + \wmem_{16} b_-^2 
\end{align*}
Then
\begin{align}
\tr(\wmem)=\wmem_1 + \wmem_2 + \wmem_3 + \wmem_4
\end{align}

The \defn{adjoint} is given by $\Ad{T}(\wmem)=[T,\wmem]$ for
$T,\wmem\in\Weil$. This vanishes if $\varepsilon=0$. However the limit
$\lim_{\varepsilon\to0} \tfrac1\varepsilon \Ad{T}(\wmem)$ is defined.
If $T$ is in the set $\{J_0,J_+,J_-,K_0\}$ then we can write
$T=\sum_{\mu,\nu=0}^{2} \chi_+^\mu T^{\mu\nu} \chi_-^\nu$. Where
$\chi^1_\pm=a_\pm$ and $\chi^2_\pm=b_\pm$.  In this case, which covers
most situations, we have
\begin{align}
\lim_{\varepsilon\to0} \tfrac1\varepsilon \Ad{T}(\wmem)
&=
\sum_{\mu,\nu=0}^{2} 
\left(
\chi_+^\mu T^{\mu\nu}\frac{\partial\wmem}{\chi_-^\nu} -
\frac{\partial\wmem}{\chi_+^\mu} T^{\mu\nu}{\chi_-^\nu}
\right) 
\label{W_Ad_eps0}
\end{align}

\begin{lemma}
\label{W_lm_trAd}
The formal trace of a basis element is given by
\begin{align}
\tr(S(s,t,u,v)) &=
2S(s-1,t-1,u,v) + 2S(s,t,u-1,v-1)
\label{W_lm_tr}
\end{align}
This commutes with the adjoint action of $J_0,J_+,J_-$ and $K_0$
\begin{align}
\tr\circ\Ad{J_0}&=\Ad{J_0}\circ\tr &
\tr\circ\Ad{J_+}&=\Ad{J_+}\circ\tr &
\tr\circ\Ad{J_-}&=\Ad{J_-}\circ\tr &
\tr\circ\Ad{K_0}&=\Ad{K_0}\circ\tr 
\label{W_lm_Ad}
\end{align}
but not with $K_+,K_-$. It also commutes with taking the conjugate 
$\tr(\wmem^\dagger)=(\tr(\wmem))^\dagger$ 
\end{lemma}

\begin{proof}
By reasoning similar to the appendix in \cite{Gratus5},
we can split the basis element for any $d\in\Intg^+$
\begin{align*}
S(u_1,u_2,u_3,u_4) = \!\!\sum_{v_1+v_2+v_3+v_4=d} \!\!
S(u_1-v_1,u_2-v_2,u_3-v_3,u_4-v_4)
S(v_1,v_2,v_3,v_4)
\end{align*}
Setting $d=2$ then 
\begin{align*}
\tr(S(s,t,u,v)) &= \wt{\tr}(S(s-1,t-1,u,v)S(1,1,0,0)) + 
\wt{\tr}(S(s,t,u-1,v-1)S(0,0,1,1)) 
\end{align*}
Since all other terms vanish hence (\ref{W_lm_tr}).

For (\ref{W_lm_Ad}) we need to set up a basis for each $\AWeil$ and
$\BWeil$ given by $S_A(s,t)$ the totally symmetric homogeneous
polynomial with each term a permutation of $a_-^s\,a_+^t$ and having
coefficient 1. Likewise for $S_B(u,v)$. These are related to
$S(s,t,u,v)$ by
\begin{align*}
(s+t)!(u+v)! S(s,t,u,v) &= (s+t+u+v)! S_A(s,t) S_B(u,v)
\end{align*}
The commutator and anti-commutator of $S_A(s,t)$ with $a_\pm$ is
given by
\begin{align*}
[a_+,S_A(s,t)] &= -\varepsilon(s+t) S_A(s,t-1) &
[a_-,S_A(s,t)] &= \varepsilon(s+t) S_A(s-1,t) \\
[a_+,S_A(s,t)]_+ &= 2s(s+t+1)^{-1} S_A(s,t-1) &
[a_-,S_A(s,t)]_+ &= 2t(s+t+1)^{-1}  S_A(s-1,t) \\
\end{align*}
Where $[\wmem_1,\wmem_2]_+=\wmem_1 \wmem_2+\wmem_2 \wmem_1$ is the
anti-commutator.  The expressions for the commutators may be proved by
induction on the degree $(s+t)$. The expressions for the
anti-commutators by use of the formula
\begin{align*}
[a_+a_-,S_A(s,t)] &= \varepsilon(s+t) S_A(s,t)
\end{align*}
Now (\ref{W_lm_Ad}) follow from direct calculation.
\end{proof}


\begin{lemma}
\label{W_lm_ab}
If $\wmem\in\Weil$ then it may be written as a sum of elements
\begin{align}
\wmem&=\sum_{r,m} \wmem_{rm}
&\text{where}&&
\Ad{K_0}\wmem &= \varepsilon r \wmem 
\,,\
\Ad{J_0}\wmem = \varepsilon m \wmem
\label{W_lm_ab_rm}
\end{align}
and we can write $\wmem_{rm}$ as a non symmetric polynomial
\begin{align}
\wmem_{rm}= \apm^{r+m} \bpm^{r-m} p_{rm}(J_0,K_0)
\end{align}
where 
\begin{align}
\apm^r = \begin{cases}
a_+^r & \text{if $r>0$} \\
1 & \text{if $r=0$} \\
a_-^{-r} & \text{if $r<0$} 
\end{cases}
\label{W_def_apm}
\end{align}
and likewise for $\bpm^r$, and where $p_{rm}(J_0,K_0)$ is a polynomial
in $J_0$ and $K_0$.
\end{lemma}

\begin{proof}
Since $\Ad{K_0}$ and $\Ad{J_0}$ commute we can diagonalise $\wmem$ with
respect to these operators. 

Since the $a_\pm$ commute with the $b_\pm$ collect all the $a_\pm$'s
together.  Use $a_-a_+=J_0+K_0+\nfrac{\varepsilon}2$ and to remove
$a_-a_+$ and $a_+a_-=J_0+K_0-\nfrac{\varepsilon}2$ to remove $a_+a_-$.
By using $[J_0,a_+]=\nfrac{\varepsilon}2a_+$ etc we can move all $J_0$
and $K_0$ to the right.  From (\ref{W_lm_ab_rm}) we know that for each
term in $\wmem$ the number of $a_+$ minus the number of $a_-$ is
$m+r$. Thus what remains on the left is $\apm^{m+r}$.  Similarly with
$b_\pm$, hence result
\end{proof}


\section{($\PsiSpace,\rho$): A Nonassociative and Noncommutative
Algebra of Formally Traceless Symmetric Polynomials in $\Weil$}
\label{ch_Pex}

In this section we introduce the subspace $\PsiSpace$ of $\Weil$ of
all formally traceless symmetric polynomials. 
We introduce two products $\rho$ and $\rho^\star$ on $\PsiSpace$
neither of which are associative, and show they are well defined. 
In lemma \ref{Pex_lm_Psi0} we show that $\PsiSpace$ contains a
subspace $\PsirSpace{0}$ which is an associative algebra, analogous to
the algebra of functions on a sphere, and that the spaces
$\PsirSpace{r}$ may be viewed as modules over $\PsirSpace{0}$. We also
show the order to which $\PsiSpace$ is associative.

\vspace{1em}

Let $\PsiSpace\subset\Weil$ be the subspace of all  traceless
symmetric polynomials in $(a_+,a_-,b_+,b_-)$.
\begin{align}
\PsiSpace=\ker(\tr)\subset\Weil
\label{Pex_def_PsiSpace}
\end{align}
We define the subspace 
\begin{align}
\PsinSpace{n}=\ker(\tr)\inter\Weil^n
\qquad\text{for $2n\in\Intg$ and $n\ge0$}
\label{Pex_def_PsinSpace}
\end{align}
i.e. the space of all formally traceless symmetric homogeneous
polynomials of order $n$. Since $\PsinSpace{n}$ is the kernel of the
restriction $\tr:\Weil^n\mapsto\Weil^{n-1}$, which is surjective, the
dimension of $\PsinSpace{n}$ is $\dim(\PsinSpace{n})=(2n+1)^2$. We
define the projection $\pi_n$
\begin{align}
\pi_n : \PsiSpace \mapsto \bigoplus_{m=0}^{2n} \PsinSpace{{m/2}} 
\label{Pex_def_pi_n}
\end{align}
We may divide the space $\PsiSpace$ in three different ways,
since for each  $\xi\in\PsiSpace$ both
$\Ad{J_0}(\xi),\Ad{K_0}(\xi)\in\PsiSpace$:
\begin{align}
\PsimSpace{m} &= 
\{\xi\in\PsiSpace\ |\ \Ad{J_0}\xi = \varepsilon m \xi\} 
\\
\PsirSpace{r} &= 
\{\xi\in\PsiSpace\ |\ \Ad{K_0}\xi = \varepsilon r \xi\}  
\end{align}
We will show in the next section that if $2n,n+m,n+r\in\Intg$ and
$n\ge0,\, |m|\le n,\, |r|\le n$ then
$\dim(\PsinSpace{n}\inter\PsirSpace{r}\inter\PsimSpace{m})=1$,
otherwise it has dimension 0.

\vskip 2em

So far we have considered $\PsiSpace$ simply as a vector space. On
this vector space we define two products. These products are given by
the projections $\rho_R$ and $\rho_R^\star$.

The intention is that $R\ge 0$ represents the radius of the sphere or
the Casimir. From (\ref{W_Casimir_J}) we want
$R^2=K_0^2-\tfrac14\varepsilon^2$.  We take the positive square root of
this equation by $K_0=\Rp=(R^2+\tfrac14\varepsilon^2)^\scrhalf$. Since
$K_0$ is not in the center of the algebra we must be careful about
quotienting by this equation. We therefore consider two projections
\begin{align}
\rho_R &: \Weil\mapsto\PsiSpace, \qquad (\rho_R)^2=\rho_R 
\notag
\\
\rho_R^\star &: \Weil\mapsto\PsiSpace, 
\qquad (\rho_R^\star)^2=\rho_R^\star
\label{Pex_def_rho_R}
\end{align}
defined by by their kernels
\begin{align}
\ker(\rho_R) &= \{\wmem\,(K_0-\Rp) \ |\
\wmem\in\Weil\} 
\notag
\\
\ker(\rho_R^\star) &= \{(K_0-\Rp)\,\wmem \ |\
\wmem\in\Weil\} 
\label{Pex_def_ker_rho_R}
\end{align}
where
\begin{align}
\Rp &= (R^2+\tfrac14\varepsilon^2)^\scrhalf, \qquad
R,\Rp\in\Real,\, R\ge0,\Rp>0
\label{Pex_def_Rp}
\end{align}
{}From now on we shall write $\rho$ and $\rho^\star$ when there is only
one possible $R$.  $\Rp$ will always be given by (\ref{Pex_def_Rp}),
Note that $\Rp\to R$ as $\varepsilon\to 0$.
These projections are surjective (onto $\PsiSpace$)
due to
\begin{lemma}
\begin{align}
\Weil = \PsiSpace \oplus \ker(\rho) = 
\PsiSpace \oplus \ker(\rho^\star)
\label{Pex_Weil_p_ker}
\end{align}
\end{lemma}

\begin{proof}
Let $\wmem\in\Weil^n$, $n\ne0$ then $\wmem K_0$ has a component in
$\Weil^{n+1}$. Furthermore this component has non zero trace, so
$\tr(\wmem K_0)$ has a component in $\Weil^n$. Now $\tr(\wmem\Rp)$
cannot have a component in $\Weil^n$ so $\tr(\wmem(K_0-\Rp))\ne 0$.
\begin{align}
\PsiSpace\inter\ker(\rho)=\{0\}
\label{Pex_lmpf_Psi_inter_ker_zero}
\end{align}
We show $\Weil = \PsiSpace \oplus \ker(\rho_R)$ by dimensional
argument.  For this proof only let 
\begin{align*}
\wh\Weil^n&=\bigoplus_{m=0}^{2n} \Weil^{m/2} &
\wh\Psi^n&=\bigoplus_{m=0}^{2n} \PsinSpace{m/2} 
\end{align*}
while we also let
\begin{align*}
\wh V^n &=\ker(\rho_R) \inter \wh\Weil^n
= \{ \wmem\,(K_0-\Rp)\ |\ \wmem\in\wh\Weil^{n-1} \}
\end{align*}
So
\begin{align*}
\dim(\wh V^n)=\dim(\wh\Weil^{n-1})=\sum_{m=0}^{2n-2}\dim(\Weil^m)
\end{align*}
Also
\begin{align*}
\dim(\wh\Psi^n) &= \sum_{m=0}^{2n}\dim\PsinSpace{{m/2}}
= \sum_{m=0}^{2n}(\dim\Weil^{m/2}-\dim\Weil^{m/2-1})
\end{align*}
Thus
\begin{align*}
\dim(\wh\Psi^n) + \dim(\wh V^n)=\dim(\wh\Weil^n)
\end{align*}
By intersecting (\ref{Pex_lmpf_Psi_inter_ker_zero}) by $\wh\Weil^n$ we
have $\wh\Psi^n \inter \wh V^n=\{0\}$.  Combining these last two
equations gives $\wh\Psi^n \oplus \wh V^n = \wh\Weil^n$.
\end{proof}

\vskip 1em

For a given $\varepsilon,R$ we define the products, also called $\rho$
and $\rho^\star$, by:
\begin{align}
\PsiSpace\times\PsiSpace &\mapsto\PsiSpace\,,
&\xi,\zeta &\mapsto \rho(\xi\zeta)\,,
&\xi,\zeta &\mapsto \rho^\star(\xi\zeta)\,,
\label{Pex_def_alg}
\end{align}
Since $\ker(\rho)$ is not a two sided ideal, the product in
$\PsiSpace$ is not associative and thus $(\PsiSpace,\rho)$ is not an
algebra.  For example
$\rho(\rho(a_+a_-)a_-)=\rho(J_0a_-+(\Rp-\nfrac{\varepsilon}{2})a_-)$
while $\rho(a_+\rho(a_-a_-))=\rho(a_+a_-a_-) =
\rho(J_0a_-+(\Rp-{\varepsilon})a_-)$

However the restriction $(\PsirSpace{0},\rho)$ is an associative
algebra. This algebra is identified with the analogue of the algebra
of functions given in \cite{Gratus5}.  Both $\PsiSpace$ and
$\PsirSpace{r}$ for all $r$ are modules overs $\PsirSpace{0}$ in the
same way the space of spinors and vector fields are modules over the
space of functions.  For reasons that will become apparent we will
call the elements of $\PsirSpace{1}$ vectors, the elements of
$\PsirSpace{-1}$ forms, and the elements of $\PsirSpace{1/2}$ and
$\PsirSpace{-1/2}$ spinors.  We also define a sesquilinear product and
quadratic map:
\begin{align}
\innerprod(\bullet,\bullet)&:\PsiSpace\times\PsiSpace\mapsto\Cmpx 
&
\innerprod(\xi,\zeta) &= \pi_0(\rho_R(\xi^\dagger \zeta))
\label{Pex_def_inner_prod}
\\
\norm{\bullet}^2 &:\PsiSpace\mapsto\Real
&
\norm{\zeta}^2 &= \innerprod(\zeta,\zeta)
\label{Pex_def_norm}
\end{align}
As we will see this sesquilinear product is hermitian but not positive
definite, and that $\norm{\zeta}^2$ may be positive negative or zero.

\begin{lemma}
\label{Pex_lm_Psi0}
For a given $\varepsilon$ and $R$, the product $\rho$ when restricted
to the space $\PsirSpace{0}$ is closed and associative
making $(\PsirSpace{0},\rho)$ an algebra. It is also equal the
product $\rho^\star$ on $\PsirSpace{0}$.

Both these algebras are equivalent to algebra $\PPexpr(\varepsilon,R)$
given by:
\begin{align}
\PPexpr(\varepsilon,R) &= \{\mbox{Free noncommuting algebra of
polynomials in $(J_+,J_-,J_0)$}\}\Big/\sim
\end{align}
where $\sim$ are the relations
\begin{align}
[J_0,J_+]&\sim \varepsilon J_+ &
[J_0,J_-]&\sim -\varepsilon J_- &
[J_+,J_-]&\sim 2\varepsilon J_0 &
J_0^2 + \tfrac12 J_+J_- + \tfrac12 J_-J_+ &= R^2 
\end{align}
This algebra is described
in detail in \textup{\cite{Gratus5}}. 

For general $\wmem_1,\wmem_2\in\PsiSpace$ we have
$\rho(\wmem_1\,\rho(\wmem_2))=\rho(\wmem_1 \wmem_2)$ whilst
$\rho(\rho(\wmem_1)\wmem_2)-\rho(\wmem_1 \wmem_2)=O(\varepsilon)$. 
In terms of the algebra $(\PsiSpace,\rho)$ this means that for
$\xi_1,\xi_2,\xi_3\in\PsiSpace$
\begin{align}
\rho(\xi_1\rho(\xi_2 \xi_3))=\rho(\xi_1\xi_2\xi_3) 
\nonumber\\
\rho(\rho(\xi_1\xi_2)\xi_3)-\rho(\xi_1\xi_2\xi_3)=O(\varepsilon)
\end{align}

If $\xi\in\PsirSpace{r_1}$ and $\zeta\in\PsirSpace{r_2}$ then
$\xi^\dagger\in\PsirSpace{-r_1}$ and
$\rho(\xi\zeta)\in\PsirSpace{r_2+r_1}$. 

Both $\PsiSpace$ and $\PsirSpace{r}$ for all $r$ are right modules
overs $\PsirSpace{0}$, since given $\xi\in\PsiSpace$ and
$f,g\in\PsirSpace{0}$ then
\begin{align}
\rho(\rho(\xi f)g) =
\rho(\xi \rho(fg)) =
\rho(\xi fg)
\end{align}
The sesquilinear product $\innerprod(\bullet,\bullet)$ is Hermitian
\begin{align}
\innerprod(\xi,\zeta) &= \cnj{\innerprod(\zeta,\xi)} \qquad \forall
\xi,\zeta\in\PsiSpace
\end{align}
\end{lemma}

\begin{proof}
If $f,g\in\PsirSpace{0}$ then it is clear that $\Ad{K_0}(fg)=0$, so
$\rho(fg)=\rho^\star(fg)$.  From lemma \ref{W_lm_ab} we can write
$fg=\apm^m\bpm^{-m}p(J_0,K_0)$, thus
$\rho(fg)=\sum_m\apm^m\bpm^{-m}p(J_0,\Rp)$. But
$\apm^m\bpm^{-m}=J_+^m$ for $m>0$ and $\apm^m\bpm^{-m}=J_-^{-m}$ for
$m<0$ so $\rho(fg)$ is a polynomial in $\{J_0,J_+,J_-\}$ and all such
polynomial are elements in $\PPexpr(\varepsilon,R)$.  From
(\ref{W_Casimir_J}) these two algebras are equivalent.

Given $\wmem_1,\wmem_2\in\Weil$ then from (\ref{Pex_Weil_p_ker}) we
may write uniquely $\wmem_1=\rho(\wmem_1)+\wmem_1'(K_0-\Rp)$. Thus
$\rho(\wmem_2\wmem_1)=
\rho(\wmem_2\rho(\wmem_1)+\wmem_2\wmem_1'(K_0-\Rp))=
\rho(\wmem_2\rho(\wmem_1))$. Whilst
$\rho(\rho(\wmem_1)\wmem_2)-\rho(\wmem_1\wmem_2)=
\rho(\wmem_1'(K_0-\Rp)\wmem_2)= \rho(\wmem_1'\Ad{K_0}(\wmem-2))=
O(\varepsilon)$

If $f,g\in\PsirSpace{0}$ and $\xi\in\PsiSpace$ then we write $\xi
f=\rho(\xi f)+ \xi_1(K_0 - \Rp)$ so $\rho(\xi fg)-\rho(\rho(\xi
f)g)=\rho(\xi_1(K_0-\Rp)g)=0$ since $[K_0,g]=0$.

Given $\xi\in\PsirSpace{r_1}$  and $\zeta\in\PsirSpace{r_2}$ then
$\innerprod(\xi,\zeta)\ne0$ only if $r_1=r_2$. In this case
$\rho^\star(\xi^\dagger \zeta)= \rho(\xi^\dagger \zeta)$, thus
$\cnj{\pi_0(\rho(\zeta^\dagger\xi))}=
(\pi_0(\rho(\zeta^\dagger\xi)))^\dagger=
\pi_0(\rho^\star(\xi^\dagger\zeta))=
\pi_0(\rho(\xi^\dagger\zeta))$

\end{proof}


\section{$\baseterm(n,r,m)$ Orthogonal basis for $\PsiSpace$}
\label{ch_Pnrm}

\begin{theorem}
\label{Pnrm_thm}
There is a natural basis of $\PsiSpace$ given by 
\begin{align}
\{ \baseterm({n},r,m) \ |\ 2n,r+n,m+n\in\Intg,\, |r|\le n, |m|\le n\}
\label{Pnrm_P_set}
\end{align}
where
\begin{align}
\baseterm({n},r,m) &= 
\varepsilon^{m-n}
\left( \frac{(n+m)!}{(2n)!(n-m)!}\right)^{\scrhalf}
(\Ad{J_-})^{n-m}(a_+^{n+r}b_-^{n-r})
\label{Pnrm_def_P}
\end{align}
This definition is extended to the case $\varepsilon=0$ by the use of
\textup{(\ref{W_Ad_eps0})}.  When written as formally traceless
symmetric homogeneous polynomials, $\baseterm(n,r,m)$ is independent
of $R$ and $\varepsilon$.  The basis elements are orthogonal with
respect to the sesquilinear form $\innerprod(\bullet,\bullet)$.
\begin{align}
\innerprod({\baseterm({n_1},r_1,m_1)},{\baseterm({n_2},r_2,m_2)})
&=
\delta_{n_1,n_2}\delta_{r_1,r_2}\delta_{m_1,m_2}
\norm{\baseterm({n_1},r_1,m_1)}^2
\label{Pnrm_innprod}
\end{align}
where $\norm{\baseterm({n},r,m)}^2$ may be positive negative or zero
and is independent of $m$.

For each $n$, the set $\{\baseterm({n},r,m)\ |\ \forall r,m\}$ form an
orthogonal basis for the set $\PsiSpace^n$.  These elements are
Eigenstates of the operators $\Ad{J_0}$, $\Ad{K_0}$, and $\Delta$.
\begin{align}
\Ad{J_0}\baseterm({n},r,m) &= \varepsilon m \baseterm({n},r,m)
\label{Pnrm_AdJ0}
\\
\Delta \baseterm({n},r,m) &= \varepsilon^2 n(n+1) \baseterm({n},r,m)
\label{Pnrm_Delta}
\\
\Ad{K_0}\baseterm({n},r,m) &= \varepsilon r \baseterm({n},r,m)
\label{Pnrm_AdK0}
\end{align}
where $\Delta = \Ad{J_0}\Ad{J_0} + \tfrac12\Ad{J_+}\Ad{J_-}
+ \tfrac12\Ad{J_-}\Ad{J_+}$,  
so
\begin{align}
\PsinSpace{n}\inter\PsirSpace{r}\inter\PsimSpace{m} &=
\spanrm\{\baseterm(n,r,m)\}
\label{Pnrm_Pnspaces}
\end{align}
The Operators $\Ad{J_+}$, and $\Ad{J_-}$ act as Ladder operators
within $\PsinSpace{n}\inter\PsirSpace{r}$
\begin{align}
\Ad{J_+}\baseterm({n},r,m) &= \varepsilon(n-m)^\scrhalf(n+m+1)^\scrhalf
\baseterm({n},r,m+1)
\label{Pnrm_AdJp}
\\
\Ad{J_-}\baseterm({n},r,m) &= \varepsilon(n+m)^\scrhalf(n-m+1)^\scrhalf
\baseterm({n},r,m-1)
\label{Pnrm_AdJm}
\end{align}
\end{theorem}

\begin{proof}
It is clear that $\baseterm(n,r,n)$ is a formally traceless symmetric
homogeneous polynomial or degree $2n$, from lemma \ref{W_lm_trAd} so
are $\baseterm(n,r,m)$ for all $m$.  Equations
(\ref{Pnrm_AdJ0}-\ref{Pnrm_AdJm}) may be proved in the same way as
the equivalent theorem in \cite[theorem 2]{Gratus5}.  Likewise for
(\ref{Pnrm_innprod}) where it is necessary to first show that
$\pi_0(\rho(\Ad{J_+}x))=\pi_0(\Ad{J_+}(\rho(\wmem)))=0$ for all
$\wmem\in\Weil$. This also comes from lemma \ref{W_lm_trAd}.
\end{proof}

A formula for the conjugate of the basis element
$(\baseterm(n,r,m))^\dagger$ is given in corollary
\ref{Hahn_herm_conj}. As stated in lemma \ref{Pex_lm_Psi0}, the
subspace $\PsirSpace{0}$ is equivalent to the algebra
$\PPexpr(\varepsilon,R)$ given in \cite{Gratus5}. By comparing the two
definitions for the basis vectors one sees that $\baseterm(n,0,m) =
P^m_n/\alpha_n$.  As a point of notation, the basis vectors here are
left unnormalised, by contrast to the notation used in the above
article. This avoids having to continually divide by a normalisation
constant, which here would depend on both $r$ and $n$.

\section{Useful Formulae for the Product $\rho$ on $\PsiSpace$}
\label{ch_form}

In this section we give a number of formulae which are useful for
various calculations and the physical interpretations. All these
formulae are extensions of the corresponding formulae for the basis of
scalars on the noncommutative sphere \cite{Gratus5}, but with added
care taken due to the nonassociative nature of $\PsiSpace$.

We start with a Hilbert space representation of $\Weil$ (subsection
\ref{ch_Hilb}).  As well as being used to calculate the value of norm
of the basis elements (subsection \ref{ch_Norm}), it is also useful
since it can be programmed into a symbolic mathematics programme such
as Maple, and hence used to calculate explicit formulae.  This is
followed by expressions to write $\baseterm(n,r,m)$ as rectangular
matrices (subsection \ref{ch_Mat}) and in terms of Hahn polynomials
(subsection \ref{ch_Hahn}).  Finally in subsection \ref{ch_Contac} we
extend \cite[Lemma 9]{Gratus5} to give a formula which is useful for
the definition of the exterior derivative.


\subsection{A Hilbert space representation of  $\Weil$}
\label{ch_Hilb}

The $\Hilb$ be a Hilbert space with basis given by $\rvec(k,j)$ where
$2k,j+k\in\Intg$ and $k\ge 0$, $|j|\le k$. The dual basis is
represented by $\lvec(k,j)$ where
$\lvec(k_1,j_1)\!\!\rvec(k_2,j_2)=\delta_{k_1k_2}\delta_{j_1j_2}$

The action of $\Weil$ on $\Hilb$ is given by
\begin{align}
a_+\rvec(k,j)&=\varepsilon^\scrhalf(k+j+1)^\scrhalf 
\rvec(k+\half,j+\half) 
&
b_+\rvec(k,j)&=\varepsilon^\scrhalf(k-j+1)^\scrhalf 
\rvec(k+\half,j-\half)
\nonumber\\
a_-\rvec(k,j)&=\varepsilon^\scrhalf(k+j)^\scrhalf 
\rvec(k-\half,j-\half) 
&
b_-\rvec(k,j)&=\varepsilon^\scrhalf(k-j)^\scrhalf 
\rvec(k-\half,j+\half) 
\label{Hil_def_action}
\end{align}
The dual action is given by
$(\lvec(k_1,j_1)\wmem)\rvec(k_2,j_2)=
\lvec(k_1,j_1)(\wmem\rvec(k_2,j_2))$
The effect of the six elements $J$s and $K$s are given by
\begin{align}
J_0\rvec(k,j) &= \varepsilon j \rvec(k,j)  
&
K_0\rvec(k,j) &= \varepsilon (k+\tfrac12) \rvec(k,j)   
\nonumber\\
J_+\rvec(k,j) &= \varepsilon(k-j)^\scrhalf(k+j+1)^\scrhalf\rvec(k,j+1)
&
K_+\rvec(k,j) &= 
\varepsilon(k-j+1)^\scrhalf(k+j+1)^\scrhalf\rvec(k+1,j)
\nonumber\\
J_-\rvec(k,j) &= \varepsilon(k+j)^\scrhalf(k-j+1)^\scrhalf\rvec(k,j-1)
&
K_-\rvec(k,j) &= \varepsilon(k-j)^\scrhalf(k+j)^\scrhalf\rvec(k-1,j)
\end{align}
One sees that the the effect of $J_+,J_-,J_0$ are exactly the same as
the elements in the algebra in \cite{Gratus5}

\begin{lemma}
\label{Hil_lm}
Given $\wmem\in\Weil$ which is independent of $R$ and $\varepsilon$,
\begin{align}
\wmem\rvec(k,j)=0\quad\forall\rvec(k,j)\in\Hilb  && \iff && 
\rho_R(\wmem)=0 \quad\forall R
&& \iff && \wmem=0
\end{align}
If $\Ad{K_0}\wmem=\varepsilon r \wmem$ and $\Ad{J_0}\wmem=\varepsilon
m \wmem$ then
\begin{align}
\wmem\rvec(k,j) \in\spanrm\{\rvec(k+r,j+m)\}
\label{Hil_lm_krjm}
\end{align}
If $R^2=\varepsilon^2 k(k+1)$ then
\begin{align}
\rho_R(\wmem)\rvec(k,j) &= \wmem\rvec(k,j)
\label{Hil_lm_rho}
\\
\pi_0(\rho_R(\wmem)) &= \frac1{2k+1}\sum_{j=-k}^k \lvec(k,j) \wmem
\rvec(k,j) 
\label{Hil_lm_pi0}
\end{align}
\end{lemma}

\begin{proof}
Given $\wmem\in\Weil$ such that $\Ad{K_0}\wmem=\varepsilon rx$ and
$\Ad{J_0}\wmem=\varepsilon mx$ we use lemma \ref{W_lm_ab} to write 
$\wmem$ in form (\ref{W_lm_ab_rm}). From (\ref{Hil_def_action}),
(\ref{Hil_lm_krjm}) is obvious.

\vskip 1em

If $\rho_R(\wmem)=0 \quad\forall R$ then it is clear that
$\wmem\rvec(k,j)=0\quad\forall j,k$. To prove the reverse
$\wmem\in\Weil$ can be written as a sum of elements
$\wmem=\sum_{rm}\wmem_{rm}$ with $\wmem_{rm}$ in the eigenspace of
$\Ad{K_0}$ and $\Ad{J_0}$. From (\ref{Hil_lm_krjm}) each $\wmem_{rm}$
satisfies $\wmem_{rm}\rvec(k,j)=0$

Using lemma \ref{W_lm_ab} write 
\begin{align*}
\wmem_{rm} &= \apm^{r+m}\bpm^{r-m}p(K_0,J_0)
\end{align*}
Then 
\begin{align*}
0 = \wmem_{rm}\rvec(k,j) &= 
p(\varepsilon(k+\tfrac12),\varepsilon j)\apm^{r+m}\bpm^{r-m}\rvec(j,k)
\qquad\forall j,k
\end{align*}
writing $\wh{p}(k,j)=p(\varepsilon(k+\tfrac12),\varepsilon j)$ then
$\wh{p}(k,j) = 0$ for all $j,k$ such that $2j,2k\in\Intg$,
$k\ge\max(0,-r)$ and $\max(-k,-k-m)\le j \le \min(k,k-m)$.  Now
writing $\wh{p}(k,j)=\sum_s p_s(j) k^s$ then $p_s(j)=0$ for above
range of $j$. Chosing $k$ greater than the degree of $\wh{p}$ implies
$p_s\equiv 0$. Thus $\wh{p}\equiv 0$, thus $\wmem=0$

\vskip 1em

Finally if $\wmem\in\PsiSpace$ writing $\wmem=\sum_{rm}\wmem_{rm}$ as
before, from (\ref{Hil_lm_krjm}) it is clear that the right hand side
of (\ref{Hil_lm_pi0}) vanishes for all but the $\wmem_{00}$ term.  Now
$\wmem_{00}$ may be written as a polynomial in $J_0$. We can use the
result of \cite[Theorem 3]{Gratus5}

\end{proof}


\subsection{The value of  $\norm{\baseterm(n,r,m)}^2$}
\label{ch_Norm}

\begin{theorem}
\label{Norm_th}
The value of $\norm{\baseterm(n,r,m)}^2$ is given by
\begin{align}
\norm{\baseterm(n,r,m)}^2
&=
\pi_0(a_-^{n+r}a_+^{n+r}b_+^{n-r}b_-^{n-r})
=
\frac{(n+r)!(n-r)!}{(2n+1)!}
\prod_{s=1}^{n-r}(2\Rp-\varepsilon s)
\prod_{s=1}^{n+r}(2\Rp+\varepsilon s)
\end{align}
The sign of this now depends on $\varepsilon$ and $R$ and is given by
\begin{align}
\mathrm{sign}(\norm{\baseterm(n,r,m)}^2)
=
\begin{cases}
1 & \text{\rm if } 2\Rp/\varepsilon>n-r
\\
0 & \text{\rm if } 2\Rp/\varepsilon\in\Intg,\ 
1\le 2\Rp/\varepsilon \le n-r
\\
(-1)^{(n-r-\lfloor 2\Rp/\varepsilon \rfloor)} &
\text{\rm if } 2\Rp/\varepsilon\not\in\Intg,\ 
0< 2\Rp/\varepsilon < n-r
\end{cases}
\end{align}
where $\lfloor 2R/\varepsilon \rfloor$ is the largest integer less
than or equal to $2R/\varepsilon$.
\end{theorem}


\begin{proof}

This proof is similar to the proof of \cite[Theorem 6]{Gratus5}.
Consider $\Rp=\varepsilon(k+\tfrac12)$, we have
\begin{align*}
\pi_0(a_-^{n+r}a_+^{n+r}b_+^{n-r}b_-^{n-r})
&=
\frac{1}{2k+1}\sum_{j=-k}^k
\lvec(k,j)a_-^{n+r}a_+^{n+r}b_+^{n-r}b_-^{n-r}\rvec(k,j)
\\
&=
\frac{\varepsilon^{2k}}{2k+1}\sum_{j=-k}^k 
\frac{(k+j+r+n)!(k-j)!}{(k+j)!(k-j-n+r)!}
\\
\intertext{Substituting $j=s-k$ and removing vanishing terms terms}
&=
\frac{\varepsilon^{2k}}{2k+1}\sum_{j=-k}^{2k-n+r} 
\frac{(s+r+n)!(2k-s)!}{s!(2k-n+r-s)!}
\\
&=
\varepsilon^{2k}\frac{(r+n)!(2k)!}{(2k+1)(2k-n+r)!}
\sum_{j=-k}^{2k-n+r}
\frac{(r+n+1)_s(-2k+n-r)_s}{s!(-2k)_s}
\\
&=
\varepsilon^{2k}\frac{(r+n)!(2k)!}{(2k+1)(2k-n+r)!}
F(r+n+1,n-r-2k;-2k;1)
\\
&=
\varepsilon^{2k}\frac{(r+n)!(2k)!}{(2k+1)(2k-n+r)!}
\lim_{\kappa\to0}
F(r+n+1,n-r-2k;\kappa-2k;1)
\\
&=
\varepsilon^{2k}\frac{(r+n)!(2k)!}{(2k+1)(2k-n+r)!}
\lim_{\kappa\to0}
\frac{\Gamma(\kappa-2k)\Gamma(\kappa-2n-1)}%
{\Gamma(\kappa-2k-r-n-1)\Gamma(\kappa-n+r)}
\\
&=
\varepsilon^{2k}\frac{(r+n)!(n-r)!(2k+1+r+n)!}{(2k+1)(2n+1)!(2k-n+r)!}
\end{align*}
Substituting $\Rp=\varepsilon(k+\tfrac12)$ and using lemma \ref{Hil_lm}
gives result.
\end{proof}


\subsection{Matrix representations for $\baseterm(n,r,m)$}
\label{ch_Mat}

If $\Rp=\varepsilon(k+\tfrac12)$ then as already noted we can quotient
the space of functions $\PsirSpace{0}$ to give the algebra of matrices
$M_{2k+1}(\Cmpx)$. This interpretation may be extended for each
$\PsirSpace{r}$ but with caution to reflect its non associative
nature. The elements of $\PsirSpace{r}$ now become rectangular
matrices with ${2k+2r+1}$ rows and $2k+1$ columns. I.e. elements of
$M_{(2k+2r+1)\times(2k+1)}(\Cmpx)$.

\begin{lemma}
Given $\Rp=\varepsilon(k+\tfrac12)$ then the mapping
\begin{align}
&\varphi^r_k : \PsirSpace{r} \mapsto M_{(2k+2r+1)\times(2k+1)}(\Cmpx) 
\nonumber\\
&(\varphi^r_k(\xi))_{\mu\nu} = 
\lvec(k+r,\mu-k-r-1) \xi 
\rvec(k,\nu-k-1)
\label{Mat_def_phi}
\end{align}
where $1\le\mu\le 2k+2r+1$ and $1\le\nu\le 2k+1$,
satisfies the following properties:
\begin{align}
&\varphi^r_k(\rho(\xi f)) =  \varphi^r_k(\xi) \varphi^0_k(f) \,,\qquad
\forall \xi\in\PsirSpace{r} \,, f\in\PsirSpace{0} 
\label{Mat_homo}
\\
&(\varphi^r_k(\xi))^\dagger \varphi^r_k(\zeta) = 
\varphi^0_k(\rho(\xi^\dagger\zeta))
\label{Mat_dagg}
\\
&\varphi^r_k(\baseterm(n,r,m))=0 \text{\upshape\qquad if } 
n \ge 2k+r+1
\label{Mat_zero}
\end{align}
\end{lemma}

\begin{proof}
Equation (\ref{Mat_homo}) follows from the use of (\ref{Hil_lm_rho})
to remove the $\rho$, simple matrix multiplication and observing that
$\sum_{\mu=1}^{2k+1} \rvec(k,\mu-k-1) \lvec(k,\mu-k-1)$ is the unit
matrix in $M_{2k+1}(\Cmpx)$. Likewise (\ref{Mat_dagg}) follows from
the same plus the recognition that
$\lvec(k,\mu-k-1)\xi^\dagger\rvec(k+r,\nu-k-r-1)=
((\varphi^r_k(\xi))^\dagger)_{\mu\nu}$

For $n\ge 2k+r+1$ then from theorem \ref{Norm_th}
\begin{align*}
0 &= \innerprod({\baseterm(n,r,m)},{\baseterm(n,r,m)}) 
= \pi_0(\rho({\baseterm(n,r,m)}^\dagger{\baseterm(n,r,m)})) \\
&= \tfrac1{2k+1}\mathrm{TR}(\varphi^0_k(\rho(\xi^\dagger\xi))) 
= \tfrac1{2k+1}\sum_{\mu=0}^{2k+2r+1} \sum_{\nu=0}^{2k+1}
|\varphi^r_k(\baseterm(n,r,m))_{\mu\nu}|^2
\end{align*}
where TR represents the matrix trace.
Hence (\ref{Mat_zero})
\end{proof}

By counting dimensions we see that the
$\ker(\varphi^r_k)=\spanrm\{\baseterm(n,r,m)\ |\ \forall m,n \text{
s.t. } n\ge 2k+r+1\}$. This is not an ideal in general but it does
have the property: if $\xi\in\ker(\varphi^r_k)$ and
$f\in\PsirSpace{0}$ then $\rho(\xi f)\in\ker(\varphi^r_k)$.


\subsection{Writing $\baseterm(n,r,m)$ in terms of Hahn polynomials}
\label{ch_Hahn}

\begin{lemma}
The highest degree term in the expansions of $\baseterm(n,r,m)$ is
given by
\begin{upshape}
\begin{align}
\baseterm(n,r,m) &= 
 (-1)^{n-\max(r,m)}
\binom{2n}{n+m}^\scrhalf
\,\apm^{r+m}\,\bpm^{r-m} 
J_0^{n-\max(|r|,|m|)} \ +\ \text{Lower Degree Terms in $J_0$}
\label{Hahn_xi_ab}
\end{align}
\end{upshape}
where $\apm$ and $\bpm$ are given by (\ref{W_def_apm})
\end{lemma}

\begin{proof}
By repeated application of the formula
\begin{align}
\Ad{J_-}(a_+^x b_+^{n+r-x} a_-^{n-r-v} b_-^v)
&=
\varepsilon x a_+^{x-1}b_+^{n-r-x+1} a_-^{n-r-v} b_-^v 
-\varepsilon v a_+^x b_+^{n+r-x} a_-^{n-r-v+1} b_-^{v-1}
\end{align}
we see that
\begin{align}
\lefteqn{(\Ad{J_-})^{n-m}(a_+^{n+r} b_+^{n-r})}\qquad
& 
\cr
&=
\varepsilon^{n-m}\sum_{s=\max(0,r-m)}^{\min(n-m,n+r)}
\frac{(n-m)!(n+r)!(n-r)!(-1)^{n-m-s}}{s!(n-m-s)!(n+r-s)!(m-r+s)!}
a_+^{n+r-s} b_+^{s} a_-^{n-m-s} b_-^{m-r+s}
\end{align}
Using $a_+a_-=J_0+\text{LDT}$ and $b_+b_-=-J_0+\text{LDT}$, where LDT
means ``Lower Degree Terms in $J_0$'', we have
\begin{align*}
a_+^{n+r-s} b_+^{s} a_-^{n-m-s} b_-^{m-r+s}
&=
\apm^{r+m}\bpm^{r-m}J_0^{n+\min(r,-m)+\min(0,m-r)}(-1)^{s+\min(0,m-r)}
+\text{LDT}
\cr
&=
\apm^{r+m}\bpm^{r-m}J_0^{n-\max(|r|,|m|)}(-1)^{s+\min(0,m-r)}
+\text{LDT}
\end{align*}
Where the last expression comes from considering each of the four case
below separately.

We have to consider separately the cases $m\ge r$ and $r\ge m$. We use
the equations $(N+s)!=N!(N+1)_s$ and $N!=(N-s)!(-1)^s(-N)_s$ to
convert our sum into a hypergeometric function. Since the coefficient
is 1 we can use the summation formula to get an explicit form in terms
of factorials.

For $m\ge r$ we have
\begin{align*}
\lefteqn{(\Ad{J_-})^{n-m}(a_+^{n+r} b_+^{n-r})}
\qquad\\
&=
\varepsilon^{n-m}
\apm^{r+m}\bpm^{r-m}J_0^{n-\max(|r|,|m|)}
\frac{(-1)^{n-m}(n-r)!}{(m-r)!}
\sum_{s=0}^{n+\min(-m,r)}
\frac{(m-n)_s(-n-r)_s}{s!(m-r+1)_s}
+\text{LDT}
\\
&=
\varepsilon^{n-m}
\apm^{r+m}\bpm^{r-m}J_0^{n-\max(|r|,|m|)}
\frac{(-1)^{n-m}(n-r)!}{(m-r)!}
F(m-n,-n-r;m-r+1;1)
+\text{LDT}
\\
&=
\varepsilon^{n-m}
\apm^{r+m}\bpm^{r-m}J_0^{n-\max(|r|,|m|)}
\frac{(-1)^{n-m}(2n)!}{(n+m)!}
+\text{LDT}
\end{align*}
whilst for $m\le r$
\begin{align*}
{(\Ad{J_-})^{n-m}(a_+^{n+r} b_+^{n-r})}
&=
\varepsilon^{n-m}
\apm^{r+m}\bpm^{r-m}J_0^{n-\max(|r|,|m|)}
\frac{(-1)^{n-r}(2n)!}{(n+m)!}
+\text{LDT}
\end{align*}
Using (\ref{Pnrm_def_P}) gives result.
\end{proof}

\begin{corol}
\label{Hahn_herm_conj}
The hermitian conjugate is given by
\begin{align}
(\baseterm({n},r,m))^\dagger &= (-1)^{(r+m)}
\baseterm({n},-r,-m)
\label{Pnrm_dagger}
\end{align}
\end{corol}

\begin{proof}
By putting $m=-n$ in above we have the exact formula 
\begin{align*}
\baseterm(n,r,-n) &= (-1)^{r+n} a_-^{n-r}b_+^{n+r}
\label{Pnrm_dagger_Pnn}
\end{align*}
Taking the conjugate of (\ref{Pnrm_def_P}), using
$(\Ad{J_-}(f))^\dagger=-\Ad{J_+}(f^\dagger)$ and substituting
(\ref{Pnrm_dagger_Pnn}) gives (\ref{Pnrm_dagger}).
\end{proof}


\begin{theorem}
We can write the polynomials $\baseterm(n,r,m)$ as
\begin{upshape}
\begin{align}
\baseterm(n,r,m) &= 
\begin{cases}
\displaystyle{
 (-1)^{n-r} \varepsilon^{n-r}
\binom{2n}{n+m}^\scrhalf\binom{2n}{n+r}^{-1}
\,a_+^{r+m}\,b_+^{r-m} \,
h^{(r-m,r+m)}_{n-r}
\left(
\frac{(J_0+\Rp)}{\varepsilon}-\frac12,\frac{2\Rp}{\varepsilon}\right)
}
\\ \text{\qquad\qquad for } -r\le m\le r 
\\
\displaystyle{
(-1)^{n-m} \varepsilon^{n-m}
\binom{2n}{n+m}^{-\scrhalf}
\,a_+^{r+m}\,b_-^{m-r} \,h^{(m-r,r+m)}_{n-m}
\left(
\frac{(J_0+\Rp)}{\varepsilon}-\frac12,\frac{2\Rp}{\varepsilon}+r-m
\right)
}
\\ \text{\qquad\qquad for }-m\le r\le m 
\\
\displaystyle{
(-1)^{n-r} \varepsilon^{n+m}
\binom{2n}{n+m}^{-\scrhalf}
\,a_-^{-r-m}\,b_+^{r-m} } \times
\\
\displaystyle{\qquad\qquad
h^{(r-m,-r-m)}_{n+m}
\left(
\frac{(J_0+\Rp)}{\varepsilon}-\frac12+r+m,\frac{2\Rp}{\varepsilon}+r+m
\right)
}
\text{\qquad for }m\le r\le -m 
\\
\displaystyle{
(-1)^{n-m} \varepsilon^{n+r}
\binom{2n}{n+m}^\scrhalf\binom{2n}{n+r}^{-1}
\,a_-^{-r-m}\,b_-^{m-r} } \times
\\
\displaystyle{\qquad\qquad
h^{(m-r,-r-m)}_{n+r}
\left(
\frac{(J_0+\Rp)}{\varepsilon}-\frac12+r+m,\frac{2\Rp}{\varepsilon}+2r
\right)
}
\text{\qquad for }r\le m\le -r 
\end{cases}
\end{align}
\end{upshape}
\end{theorem}

\begin{proof}
It is clear we shall have to consider the four cases separately. For
each case the proof is not too different from \cite[Theorem
7]{Gratus5}.
Let us here prove the theorem for the first case where $-r\le m\le
r$, the other cases are Similarly proved. 

Given two basis elements $\baseterm(n_1,r_1,m_1)$ and
$\baseterm(n_2,r_2,m_2)$ it is clear that
$\innerprod({\baseterm(n_1,r_1,m_1)},{\baseterm(n_2,r_2,m_2)})=0$
unless $m_1=m_2$ and $r_1=r_2$.  Use lemma (\ref{W_lm_ab}) so the
$\baseterm(n_1,r,m)=a_+^{r+m} b_+^{r-m} p(J_0,K_0)$. Thus
$\baseterm(n_1,r,m)=\rho(\baseterm(n_1,r,m))=a_+^{r+m} b_+^{r-m}
p_{n_1}(J_0)$ for the polynomial $p_{n_1}(J_0)$ which we want to
determine.

Let $R^2=\varepsilon^2 k(k^1)$ so $R'=\varepsilon(k+\tfrac12)$
Then using (\ref{Hil_lm_pi0}) we have for $n_1\ne n_2$
\begin{align*}
0 &= \innerprod({\baseterm(n_1,r,m)},{\baseterm(n_2,r,m)})
\\
&=
\sum_{j=-k}^{k}
\lvec(k,j) 
p_{n_1}(J_0) a_+^{r+m} b_+^{r-m} 
a_+^{r+m} b_+^{r-m} p_{n_2}(J_0)
\rvec(k,j)
\\
&=
\sum_{j=-k}^{k}
p_{n_1}(\varepsilon j) 
p_{n_2}(\varepsilon j) 
\frac{(k+j+r+m)!(k-j+r-m)!}{(k+j)!(k-j)!}
\end{align*}
Thus $p_{n}(\varepsilon j)$ are polynomials orthogonal with respect to
a weight function. This weight function is precisely the same as that
for the Hahn polynomials $h^{(\alpha,\beta)}_{n'}(x,N)$ where
$\alpha=r-m$, $\beta=r+m$, $x=j+k$ and $N=2k+1$. We see that $n'=n-r$
since $h^{(\alpha,\beta)}_{n'}(x,N)$ is a polynomial of order $n'$ and
in $x$, and from (\ref{Hahn_xi_ab}) $p_n(J_0)$ is a polynomial of order
$n-r$.  Comparing the coefficient for the $J_0^{n-r}$ term we have
{}from \cite{Nik}
\begin{align*}
h^{(r-m,r+m)}_{n-r}
\left(\frac{(J_0+\Rp)}{\varepsilon}-\frac12,\frac{2\Rp}{\varepsilon}
\right)
&=
\frac{(2n)!}{(n-r)!(r+n)!}\frac{J_0^{n-r}}{\varepsilon^{n-r}}
+\text{LDT}
\end{align*}
Hence result for the range $-r\le m\le r$. Since this is true for all
$k$ then using lemma \ref{Hil_lm} gives result for all $\varepsilon,R$.
The other ranges are proved similarly.
\end{proof}

In the limit $\varepsilon=0$ the Hahn polynomials become Jacobi
Polynomials $(P^{(\alpha,\beta)}_n(z))$. Thus we have
\begin{align}
\baseterm(n,r,m)|_{\varepsilon=0} = 
(-1)^{n-\max(r,m)} (2R)^{n-\max(|r|,|m|)}
\binom{2n}{n+m}^\scrhalf\!&\binom{2n}{\!n\!-\!\max(|r|,|m|)\!}^{-1}
\times
\nonumber\\ &
\apm^{r+m}\bpm^{r-m} \, P^{(|r-m|,|r+m|)}_{n-\max(|r|,|m|)}
\left(\frac{J_0}{R}\right)
\label{Hahn_eps0}
\end{align}


\subsection{Another useful theorem}
\label{ch_Contac}

Here we extend \cite[Lemma 9]{Gratus5} to give a formula which is
useful for the definition of the exterior derivative.

\begin{theorem}
\label{th_Contac}
For $r_1,r_2,n$ let
\begin{align}
&\omega^{r_1r_2}_n:\PsiSpace\mapsto\PsiSpace
\cr
&\omega^{r_1r_2}_n(\xi)= \sum_{m=-n}^n 
\rho(\baseterm(n,r_2,m)^\dagger \, \xi \,\baseterm(n,r_1,m))
\end{align}
Then $\omega^{r_1r_2}_n$ commutes with 
$\Ad{J_0}$, $\Ad{J_+}$, $\Ad{J_-}$, $\Delta$ and satisfies 
\begin{align}
\Ad{K_0}\omega^{r_1r_2}_n(\xi)=\varepsilon(r_1-r_2)
\omega^{r_1r_2}_n(\xi) +
\omega^{r_1r_2}_n(\Ad{K_0}\xi)
\end{align}
It effect on the basis elements is given by
\begin{align}
\omega^{r_1r_2}_{n_1}(\baseterm(n,r,m))
&=
\omega^{r_1r_2r}_{n_1 n}\baseterm(n,r+r_1-r_2,m)
\end{align}
where $\omega^{r_1r_2r}_{n_1 n}\in\Cmpx$. For $\xi=1$
\begin{align}
\omega^{r_1r_2}_n(1)=\delta_{r_1r_2} (2n+1)
\norm{\baseterm(n,r_1,m)}^2
\end{align}
\end{theorem}

\begin{proof}
The effect of $\Ad{J_0}$, $\Ad{J_+}$, $\Ad{J_-}$, and $\Ad{K_0}$ come
{}from direct calculation using (\ref{Pnrm_AdJ0}), (\ref{Pnrm_AdJp}),
(\ref{Pnrm_AdJm}) and  (\ref{Pnrm_AdK0}) respectively. The effect of
$\Delta$ comes from $\Ad{J_+}$ and $\Ad{J_-}$.  
{}From this, it is clear that $\omega^{r_1r_2}_{n_1}(\baseterm(n,r,m))$
must be proportional to $\baseterm(n,r+r_1-r_2,m)$ and that the
proportionality constant must be independent of $m$ since
$\omega^{r_1r_2}_{n_1}$ commutes with $\Ad{J_+}$. 

Since $\omega^{r_1r_2}_n(1)$ is proportional to
$\baseterm(0,r_1-r_2,0)$ it is clear that it vanishes unless
$r_1=r_2$. In this case taking $\pi_0$ in the expression give the
result.
\end{proof}


\section{Physical Interpretation}
\label{ch_Phys}

There are two ways of interpreting the algebra $(\PsiSpace,\rho)$ and
the basis elements $\baseterm(n,r,m)$. In subsection \ref{ch_drm} we
show that $\baseterm(n,r,m)$ my be thought of as the deformed rotation
matrices for $SU(2)$, whilst in subsections \ref{ch_covec},
\ref{ch_vec} and \ref{ch_spin} we use them to construct the analogue
of scalar, vector and spinor fields.


\subsection{Deformed rotation matrices}
\label{ch_drm}

In the limit $\varepsilon=0$ the elements of $\PsiSpace$ may be
regarded as complex functions over $S^3$ or $SU(2)$. The product
$\rho$ becomes pointwise multiplication.  The basis elements
$\baseterm(n,r,m)$ are proportional to the rotation matrix $D^n_{mr}$.
Thus we may interpret $\baseterm(n,r,m)$ as the deformed (or
noncommutative) rotation matrices.

Just as in standard quantum mechanics the non commutativity in the
quantum domain becomes the Poisson Bracket in the commutative case we
can define a bracket (\ref{drm_def_brak}) as the limit of the
commutator. However, since for $\varepsilon\ne0$, $\rho$ is not
associative, this bracket is neither a derivation nor satisfies the
Jacobi identity. Therefore it is not a Poison Bracket. This
is reflected in the exisitance of the number operator $\delta_N$ in
(\ref{drm_res_brak}). It would be nice if one could find a physical
interpretation for this bracket.

\begin{lemma}
\label{drm_lm}
In the limit $\varepsilon=0$ the noncommutative, nonassociative algebra
$(\PsiSpace,\rho)$ becomes the commutative associative algebra
$C_0(SU(2))$. This is the algebra of complex valued polynomial
functions on $SU(2)$ with pointwise multiplication.  

If we use the Euler angles
$(\alpha,\beta,\gamma)$ to parameterise $SU(2)$ then the generators of
$\PsiSpace$ maybe written
\begin{align}
a_+ &= \sqrt{2R} \cos(\beta/2) e^{-i/2(\alpha+\gamma)}
&
b_+ &= -i\sqrt{2R} \sin(\beta/2) e^{i/2(\alpha-\gamma)}
\nonumber\\
a_- &= \sqrt{2R} \cos(\beta/2) e^{i/2(\alpha+\gamma)}
&
b_- &= i\sqrt{2R} \sin(\beta/2) e^{-i/2(\alpha-\gamma)}
\label{drm_ab}
\end{align}
The basis elements of $\PsiSpace$ are given by the rotation
matrices
\begin{align}
\baseterm(n,r,m)|_{\varepsilon=0} &=
(i)^{m-r}(-1)^{n-r}(2R)^n \binom{2n}{r}^{-\scrhalf}
D^n_{mr}(\alpha,\beta,\gamma)
\label{drm_base}
\end{align}
where $D^n_{mr}(\alpha,\beta,\gamma)$ uses the notation in
\textup{\cite{Biedenharn1}}. The bilinear form on $\PsiSpace$ becomes
a true inner product and is given by
\begin{align}
\innerprod(\xi,\zeta) &=
\frac1{8\pi^2}\int_0^{2\pi} d\alpha \int_0^{2\pi} d\gamma 
\int_0^{\pi} d\beta
\sin(\beta)\, \cnj{\xi}\zeta 
\label{drm_innprod}
\end{align}

There is a bracket defined by
\begin{align}
&\{\bullet,\bullet\} : C_0(SU(2)) \times C_0(SU(2)) \mapsto C_0(SU(2)),
&
&\{\xi,\zeta\}=\lim_{\varepsilon\to0}\tfrac1{i\varepsilon} \rho
([\xi,\zeta])
\label{drm_def_brak}
\end{align}
which is given by
\begin{align}
\{\xi,\zeta\} &=
\frac{1}{R\sin\beta}\left( 
\frac{\partial \xi}{\partial \alpha}
\frac{\partial \zeta}{\partial \beta} 
-
\frac{\partial \xi}{\partial \beta}
\frac{\partial \zeta}{\partial \alpha} 
\right) +
\left( 
{\delta_N (\xi)}
\frac{\partial \zeta}{\partial \gamma} 
-
\frac{\partial \xi}{\partial \gamma}
{\delta_N (\zeta)} 
\right) +
\frac{\cot\beta}{R}\left( 
\frac{\partial \xi}{\partial \beta}
\frac{\partial \zeta}{\partial \gamma} 
-
\frac{\partial \xi}{\partial \gamma}
\frac{\partial \zeta}{\partial \beta} 
\right) 
\label{drm_res_brak}
\end{align}
where $\delta_N$ is the number operator
\begin{align}
&\delta_N : C_0(SU(2)) \mapsto C_0(SU(2))\, , 
&&\delta_N (\baseterm(n,r,m))= n\baseterm(n,r,m)
\label{drm_def_delN}
\end{align}

\end{lemma}

\begin{proof}
By setting $\varepsilon=0$ it is clear that the algebra is commutative
and associative. The condition that $K_0=R$ implies that the elements
of $\PsiSpace$ are functions on $S^3$, or alternatively $SU(2)$.
Assuming (\ref{drm_ab}) then (\ref{drm_base}) follows from
(\ref{Hahn_eps0}). The inner product (\ref{drm_innprod}) follows
{}from the orthogonality of the rotation matrix elements.

In the limit $\varepsilon=0$ the algebra $\Weil(0)$ is the algebra
$C_0(\Real^4)$ of complex valued polynomial functions over $\Real^4$.
The Posion Bracket corresponding to the limit of the commutator is
given by 
\begin{align*}
\{f,g\}_\Weil &= 
\frac{\partial f}{\partial a_+}
\frac{\partial g}{\partial a_-}
-
\frac{\partial f}{\partial a_-}
\frac{\partial g}{\partial a_+}
+
\frac{\partial f}{\partial b_+}
\frac{\partial g}{\partial b_-}
-
\frac{\partial f}{\partial b_-}
\frac{\partial g}{\partial b_+}
\end{align*}
Writing this with respect to the coordinates $(R,\alpha,\beta,\gamma)$
gives (\ref{drm_res_brak}) but with $\delta_N$ replaced by
$\partial/\partial R$. However, for every polynomial
$\wmem(a_+,a_-,b_+,b_-)$ not explicitly dependent on $R$, $\partial
\wmem/\partial R= \delta_N(\wmem)$. Thus restriction of
$\{\bullet,\bullet\}_\Weil$ to $C_0(SU(2))$ is valid.
\end{proof}


\subsection{The exterior derivative and 1-forms}
\label{ch_covec}

In standard geometry there are many equivalent definitions of a vector
field. However not all of the can be extended to the noncommutative
case at the same time. If ones take the definition of a vector field
given that it must satisfies the Leibniz formula, and extends this
definition to noncommutative geometry then then, at least for the
matrix case, the vector fields are given by $X=\Ad{f}$ for some
$f\in\PsirSpace{0}$. This can then be used to give a definition of the
exterior derivative and the exterior algebra \cite{Madore_book}.

Here we give another definition of covectors by use of the analogue of
the follow definition of the exterior derivative of scalar fields:
\begin{align}
df &= \sum_i \frac{\partial f}{\partial x^i} dx^i
\label{Vec_df_partial}
\end{align}
If we consider the noncommutative sphere as a three dimensional
manifold with normalised basis coordinates
$x^m=(2\Rp+\varepsilon)^{-\scrhalf}\baseterm(1,0,m)$ for
$m\in\{-1,0,1\}$ we can define the basis 1-forms as
\begin{align}
dx^m=(2\Rp+\varepsilon)^{-\scrhalf}\baseterm(1,-1,m)
\qquad\text{for } m\in\{-1,0,1\}
\label{Vec_def_dx}
\end{align}
This enables us to define the \defn{exterior derivative} on the space
of functions as
\begin{align}
&d:\PsirSpace{0}\mapsto\PsirSpace{-1}
&
df &= (2\Rp+\varepsilon)^{-1}\omega^{01}_1(f)
\qquad\forall f\in\PsirSpace{0} 
\label{Vec_def_df}
\end{align}
Since $\omega^{01}_1(1)=0$, we see that
\begin{align}
df &= \sum_{m=-1}^1 (-1)^{m+1} dx^{-m} \Ad{x^m} f
\end{align}
which is analogous to (\ref{Vec_df_partial}). This obeys the
following:

\begin{lemma}
The two definition of the 1-forms $dx^m$
given by (\ref{Vec_def_dx}) and (\ref{Vec_def_df}) agree.
For the basis elements which are functions
\begin{align}
d\baseterm(n,0,m)&=2
n(\Rp+\varepsilon n/2)(\Rp+\varepsilon)^{-1}\baseterm(n,-1,m)
\label{Vec_d_xi_nm}
\end{align}
And in the limit as $\varepsilon\to0$
\begin{align}
d(fg)=d(f)g + f (dg) + O(\varepsilon)
\label{Vec_d_limit}
\end{align}
\end{lemma}

\begin{proof}
{}From there definition we see that
\begin{align*}
J_+ &= \baseterm(1,0,1) &
-\sqrt{2}J_0 &= \baseterm(1,0,0) &
-J_- &= \baseterm(1,0,-1) 
\\
b_-^2 &= \baseterm(1,-1,1) &
\sqrt{2}a_-b_- &= \baseterm(1,-1,0) &
a_-^2 &= \baseterm(1,-1,-1) 
\end{align*}
By manipulations of $J_+,J_-$ we can show that 
\begin{align*}
[J_-,J_+^n]=-2\varepsilon J_+^{n-1}(nz+\varepsilon n(n-1)/2)
\end{align*}
{}From theorem \ref{th_Contac} we know that
$\omega^{01}_1(\baseterm(n,-1,m))$ is proportional to
$\baseterm(n,0,m)$ and its coefficient is independent of $m$. Therefore
taking $m=n\ge 1$ we have
\begin{align*}
\omega^{01}_1(\baseterm(n,0,n)) &=
a_-^2\Ad{J_+}(J_+^n) +
2a_-b_-\Ad{J_0}(J_+^n) +
b_-^2\Ad{J_-}(J_+^n) 
\\
&=
2n\varepsilon a_-b_-J_+^n -
2\varepsilon b_-^2 J_+^{n-1}(nz+\varepsilon n(n-1)/2)
\end{align*}
Now $b_-^2J_+^{n-1}=a_+^{n-1}b_-^{n+1}=\baseterm(n,-1,n)$
whilst
\begin{align*}
a_-b_-a_+^nb_-^n &=
(J_0+K_0+\varepsilon/2) a_+^{n-1}b_-^{n+1}
\\
&=
a_+^{n-1}b_-^{n+1} (J_0+K_0+n\varepsilon-\varepsilon/2))
\end{align*}
Hence (\ref{Vec_d_xi_nm}).  We show that the two definitions of $d x^m$
are equivalent by substituting $n=1$ into this formula. Finally we
note that
\begin{align*}
\sum_{m=-1}^1 (
\baseterm(1,1,m)^\dagger fg \baseterm(1,0,m)
-
\baseterm(1,1,m)^\dagger f \baseterm(1,0,m) g
-
f \baseterm(1,1,m)^\dagger g \baseterm(1,0,m) )
&=
[\baseterm(1,1,m)^\dagger ,f][g ,\baseterm(1,0,m)]
= O(\varepsilon^2)
\end{align*}
\end{proof}


\subsection{Vector fields}
\label{ch_vec}

We can now define the vectors fields as the set $\PsirSpace{1}$. The
action of a vector on a scalar is given by
\begin{align}
X (f) &= \rho((df)\,X)
\end{align}
By the same reasoning as (\ref{Vec_d_limit}), this is also only a
derivation in the limit.
\begin{align}
X(fg) &= X(f) g + f X(g) + O(\varepsilon)
\end{align}
However if we define $X_m$ to be the conjugate of $dx^m$,
$X_m=(dx^m)^\dagger=
(-1)^{m+1}(2\Rp+\varepsilon)^{-\scrhalf}\baseterm(1,1,-m)$ then 
{}from $\omega^{1,0}_1(1)=\omega^{0,-1}_1(1)=0$
we have 
\begin{align}
\sum_{m=-1}^1 x^m X_m = \sum_{m=-1}^1 X_m x^m =  0
\end{align}
identically for all $\varepsilon$ and not just in the limit
$\varepsilon=0$.  

We can now define a \defn{metric}
\begin{align}
&g:\PsirSpace{1}\times\PsirSpace{1}\mapsto\PsirSpace{0}
&& g(X,Y) = \rho(X^\dagger Y)
\end{align}
This can be extended for all spaces $\PsirSpace{r}$. 

\subsection{Spinor fields}
\label{ch_spin}

It is natural now to define the set $\PsirSpace{1/2}$ as the space of
spinor fields and its dual $\PsirSpace{-1/2}$. This means that a
vector is the product of two spinors.  We may also consider the space
$\PsinrSpace{n}{r}=\PsinSpace{n}\inter\PsirSpace{r}$ as a $2n+1$
dimensional representation of $su(2)$ given by the adjoint
representation, $\Ad{u}:\PsinrSpace{n}{r}\mapsto\PsinrSpace{n}{r}$ for
any non zero element $u\in su(2)$. The $2n+1$ eigenvalues of this
mapping are given by the set $\{m\norm{u}/2\ |\ m+n\in\Intg,|m|\le n$,
where $\norm{u}$ corresponds to the killing form. 
Thus
\begin{align*}
e^{\pi i u/\norm{u}} \baseterm(n,r,m) e^{-\pi i u/\norm{u}} =
(-1)^{2n} \baseterm(n,r,m) 
\end{align*}
Hence rotation by $2\pi$ does not change the sign of $\baseterm(n,r,m)$
if $n$ is an integer. i.e. for scalars, vectors and other ``Bosons''.
Whilst $\baseterm(n,r,m)$ changes sign under rotation of $2\pi$ for
spinors tensor fields and other ``Fermions''.

An alternative way to define spinor fields 
is as the set
\begin{align}
\spinset &= 
\left\{\binom{a_+}{b_+}f_1 + \binom{a_-}{b_-}f_2 \ \bigg|\
f_1,f_2\in\PsirSpace{0}\right\}
\label{Spin_def_mat}
\end{align}
is usually regarded as the space of spinors. This is decomposed into
$\spinset=\spinset_+\oplus\spinset_-$ corresponding to the eigenspaces
of $\Ad{K_0}$ which is now regarded as the Chirality operator.  This
interpretation is now equivalent to the one proposed by Grosse et al
\cite{Grosse5,Grosse3,Grosse6}, who go on to define and solve the
Dirac equation.  We note that $\spinset$ is not equivalent to simply
two copies of $\PsirSpace{-1/2}\oplus\PsirSpace{1/2}$, but a proper
subset. This is because the element
$\displaystyle{\binom{a_+}{0}\not\in\spinset}$.  In the commutative
limit $\varepsilon=0$ we can replace
$a_+=a_-=R^\scrhalf\cos(\theta/2)$ and
$b_+=\cnj{b_-}=R^\scrhalf\sin(\theta/2)e^{-i\phi}$ to obtain the
standard spinors on a sphere of radius $R$ and coordinates
$(0\le\theta\le\pi,0\le\phi\le 2\pi)$.


\section{Problems and Outlook}
\label{ch_Prob}

There are many problems with this interpretation of our system. Here
is a list of cases where this noncommutative geometry is different from
the commutative geometry even in the limit $\varepsilon=0$:

(1) the space $\PsirSpace{-1}$ is the image of $\PsirSpace{0}$ under
$d$. This means that all 1-forms are closed. 

(2) The definition (\ref{Vec_def_df}) can be extended as a map 
$d:\PsirSpace{r}\mapsto\PsirSpace{r-1}$ for all $r$. However this map
does not satisfy $d^2=0$ even in the limit. As a result no exterior
calculus is defined here.

(3) The set of vectors $X_i$ are not dual, in the usual sense, to the
set of 1-forms $d x^j$, or alternatively the vectors $X_i$ are not
orthogonal with respect to the metric $g$. This is because for $i\ne j$
\begin{align*}
g(X_{i},X_{j}) = \rho(dx^{i} X_{j}) = X_{j}(x^{i}) \ne 0 
\end{align*}
for all $\varepsilon$ even when $\varepsilon=0$. However we do have
$\pi_0(g(X_{i},X_{j}))=\delta_{ij}$ 

Some of these problems, together with the non derivative nature of $d$
may be solved by redefining the space of covectors. For example it may
be similar to (\ref{Spin_def_mat}).

As stated at the end of subsection \ref{ch_drm} the
bracket (\ref{drm_def_brak}) is in need of an interpretation.

\vskip 1em

Since $\rho$ is not associative it would be nice to know what rules
(if any) it obeys. Note, for example even $(\rho(\xi\xi)\xi)\ne
\rho(\xi\xi\xi)$ so $\xi^n$ is not well defined.

In the limit $\varepsilon=0$ there is a formula for the combination of
two matrix rotation entries $D^n_{rm}$. In the general case we can use
the Wigner-Eckart theorem to give
\begin{align*}
\baseterm(n_1,r_1,m_1) \baseterm(n_2,r_2,m_2)
&= \sum_{n=|n_1-n_2|}^{n_1+n_2} 
C^{n_1,n_2,n}_{m_1,m_2,m_1+m_2} 
R^{n_1,n_2,n}_{r_1,r_2,r_1+r_2} \baseterm(n,r_1+r_2,m_1+m_2)
\end{align*}
where $C^{n_1,n_2,n}_{m_1,m_2,m_1+m_2}$ is the Wigner or
Clebsch-Gordan coefficient and $R^{n_1,n_2,n}_{r_1,r_2,r_1+r_2}$ the
corresponding reduced matrix element. It would be useful to have a
formula for such elements. We note, however, that one cannot use the
Wigner-Eckart theorem to reduce the problem on the $r$ index since
this is not raised and lowered by an action of $su(2)$.

\vskip 1em

Since the space of spinor and vector fields naturally form  modules
over the space of scalar fields one can use this as a stating point
for algebraic connections \cite{Coquereaux1}.

As stated $\Weil$ contains a representation of $su(1,1)$. If one
quotiented $\Weil$ by the ideals created by $J_0\sim \Rp$, one would
gain a very similar structure which could be interpreted as the basis
of scalar, spinor and vector fields on the noncommutative analogue of
either the two dimensional DeSitter space or the Hyperbolic plane.

The Jordan-Schwinger representation can be used for any Lie algebra.
Hence one should be able to construct spinor fields on any symmetric
manifold. The combination of these two steps would give one a
description of scalar, spinor and vector fields on more exotic
symmetric spaces, such as the Einstein-DeSitter Universe.


\subsection*{Acknowledgements}

The author would like to thank the Royal Society of London for a
European Junior Fellowship which enabled him to study in at Paris.
The author would also like to thank Peter Pre\v{s}najder, John Madore,
Luiz Saeger and Jihad Mourad for useful discussions which aided
this work and, especially, Richard Kerner and the Laboratoire de
Gravitation et Cosmologie Relativistes, Paris~VI for their
hospitality.



\end{document}